\documentclass[preprint,12pt]{elsarticle}

\usepackage{amssymb}
\usepackage{amsmath}
\usepackage{amsthm}
\newtheorem{lemma}{Lemma}
\newtheorem{proposition}{Proposition}
\journal{Journal of Multivariate Analysis}

\begin{document}

\begin{frontmatter}



\title{Probability Distributions for Counts and Compositions}


\author{Guanyi Wu} %% Author name

\affiliation{organization={University of Michigan, Department of Biostatistics},%Department and Organization
            addressline={}, 
            city={Ann Arbor},
            postcode={48104}, 
            state={Michigan},
            country={United States}}

%% Abstract
\begin{abstract}
This article present a method of mutual transformation between count model and composition model. Offer the mathematical view of classical radio and log-radio in compositional data analysis and expand the idea of mixture model of counts data to the case of compositional data.
\end{abstract}




%% Keywords

\end{frontmatter}



\section{Introduction}
\label{sec1}
Modeling non-negative counts and compositional data require probability models that respect discreteness, sparsity, over-dispersion and constant-sum constrains\cite{zhang2024}. The most used Poisson and Multinomial models are often too rigid in real data\cite{holmes2012}; mixture and transformation, such as Poisson-Gamma mixture and alr for compositions, expand variability while preserving interpretability.

Specifically, for compositional data, we revisit the Dirichlet family and show what ratio and log-radio transforms of a Dirichlet random vector will yield. And we discuss the coordinate on simplex, under the view of manifold structure and chart, to justify dropping the redundant coordinate and compute Jacobian, which is connecting to Aichison's framework and the isometric log-radio(ILR) transform\cite{egozcue03}.

And for transforming counts into compositions, we use a "normalize-condition" method that maps a count vector to the simplex and expand the mixture model of counts to compositions. Let $S=\sum_{i=1}^nX_i$ is the total counts, we first model on $S$ which is Poisson-Gamma mixture, equivalent to Negative Binomial distribution. Then we fix the total count and establish the distribution of $X$ under fixed counts and intensity is Multinomial. Next we integrate the latent parameter which is Dirichlet and get $X$ conditional on $S$ is Dirichlet Multinomial, whose marginal is Beta-Binomial. This process maps mixture model of counts to compositions.

For applications, microbiome studies offer use cases. When multinomial counts are overdispersed, Dirichlet–multinomial (DM) models and their regression variants better capture between-sample variability; within this framework, Chen, Li proposed a sparse group-L1 DM regression with LRT for high-dimensional covariate selection at the taxon level, showing superior performance to models that ignore overdispersion or use proportions only\cite{chen2013}.

\section{Dirichlet Distribution with ratio}
Dirichlet distribution $X=(X_1, \ldots, X_n) \mathrel{\sim} \mathrm{Dir}(\alpha)$ is a kind of multivariate distribution on simplex, typically regarded as normalized sum of gamma distribution\cite{townes20}, with all components are between $(0,1)$ and the sum is equal to $1$.The joint PDF of X is given by
$$f_{X}(x)=\frac{1}{B(\alpha)}\prod_{i=1}^{n}x_i^{\alpha_{i}-1},$$
where $B(\alpha)=\frac{\prod_{i=1}^{n}\Gamma(\alpha_i)}{\Gamma(\sum_{i=1}^{n}\alpha_i)}$.

Now we are interested in the joint PDF of ratio random variable $Y=(\frac{X_1}{X_n},\ldots,\frac{X_{n-1}}{X_n})$, which can be seen as a transformation of $X$. Consider the multivariate transformation $Y=g(X)$. This is an invertible transformation and the inverse function is given by $X_i=\frac{Yi}{1+\sum_{i=1}^{n-1}Y_i}$, $i=1\ldots n-1$ and $X_n=\frac{1}{1+\sum_{i=1}^{n-1}Y_i}$. For convenience, we will denote $Z=1+\sum_{i=1}^{n-1}Y_i$ in the rest of the paper. Then we can calculate jacobian of the inverse function and get PDF of $Y$ by Theorem of change of variable. But the case here is somewhat different, the derivative respect to each $Y_i$ will yield a matrix with rank $n-1$. At this point, many references simply state that the last component $\frac{1}{1+\sum_{i=1}^{n-1}Y_i}$ is redundant since the restriction on simplex, like George and Irwin \cite{tiao1965}. Here we demonstrate a detailed derivation for clarity. Consider $S^{n-1}:=\{(X_1,\ldots,X_n)|\sum_{i=1}^{n-1}X_i<1\}$ and $\mathbb{R}_{+}^{n-1}$ are supports of $X$ and $Y$ respectively, and both are smooth manifold\cite{egozcue03}. Then the inverse function $g^{-1}$ can be seen as smooth map between $S^{n-1}$ and $\mathbb{R}_{+}^{n-1}$. Our aim is to obtain Jacobi matrix of this map, the first step is to equip manifolds with proper charts. Consider $T_pS^{n-1}$ and $T_{g^{-1}(p)}\mathbb{R}^{n-1}_{+}$ are tangent space of $S^{n-1}$ and $\mathbb{R}^{n-1}_{+}$ at point $p$ of $S^{n-1}$. Then $$dg^{-1}_{p}:T_pS^{n-1}\longrightarrow T_{g^{-1}(p)}\mathbb{R}^{n-1}_{+}$$ is a differential of $g^{-1}$. Take $$\phi(X_1,\ldots,X_n)=(X_1,\ldots,X_{n-1})\quad and\quad \psi(Y_1,\ldots,Y_{n-1})=(Y_1,\ldots,Y_{n-1})$$ as charts of $S^{n-1}$ and $\mathbb{R}^{n-1}_{+}$, we can then establish the Jacobi of $g^{-1}$:

\begin{align*}
    J_{g^{-1}}&=\frac{1}{z^2}\begin{pmatrix}
z-y_1&-y_1&\ldots&-y_1\\
-y_2&z-y_2&\ldots&\vdots\\
\vdots&\vdots&\ddots&\vdots\\
-y_{n-1}&-y_{n-1}&\ldots&z-y_{n-1}
\end{pmatrix}\\
&=\frac{1}{z^2}(zI_{n-1}-\begin{pmatrix}
    y_1\\
    \vdots\\
    y_{n-1}
\end{pmatrix}
\begin{pmatrix}
    1&\ldots&1
\end{pmatrix})
\end{align*}

For matrix with the form above, we can apply an useful Lemma to obtain its determinant.
\begin{lemma}
    Suppose $A$ is an invertible square matrix and $u$, $v$ are column vectors. Then the matrix determinant lemma states that 
    $$\mathrm{det}(A+uv^T)=(1+v^TA^{-1}u)\mathrm{det}(A).$$
\end{lemma}

Denote $A=\frac{1}{z}I_{n-1}$, $u=-\frac{1}{z^2}\begin{pmatrix}
    y_1\\
    \vdots\\
    y_{n-1}
\end{pmatrix}$ and $v^T=\begin{pmatrix}
    1&\ldots&1
\end{pmatrix}$, then 
\begin{align*}
\mathrm{det}(J_{g^{-1}})&=(1-\frac{1}{z^2}\begin{pmatrix}
    1&\ldots&1
\end{pmatrix}\begin{pmatrix}
    z&0&\ldots&0\\
    0&z&\ldots&0\\
    \vdots&\vdots&\ddots&\vdots\\
    0&0&\ldots&z
\end{pmatrix}\begin{pmatrix}
    y_1\\
    \vdots\\
    y_{n-1}
\end{pmatrix})(\frac{1}{z})^{n-1}\\
&=(1-\frac{\sum^{n-1}_{i=1}y_i}{z})(\frac{1}{z})^{n-1}\\
&=(\frac{1}{z})^n
\end{align*}

Thus, the PDF of $Y$ is 
\begin{align*}
    f_Y(y)&=f_X(g^{-1}(y))*|J_{g^{-1}}|\\
    &=\frac{1}{B(\alpha)}\prod_{i=1}^{n-1}(\frac{y_i}{z})^{\alpha_i-1}\cdot(\frac{1}{z})^{\alpha_n-1}\cdot(\frac{1}{z})^n\\
    &=\frac{1}{B(\alpha)}\prod_{i=1}^{n-1}y_i^{\alpha_i-1}\times z^{-\sum_{i}^n\alpha_i},
\end{align*}
which is Inverted Dirichlet Distribution.

\section{Dirichlet Distribution with log-ratio}
In this section we still assume $X\mathrel{\sim}\mathrm{Dir}(\alpha)$, while give insight on the case of log-ratio, that is, $Y=h(X)=(\log(\frac{x_1}{x_n}),\ldots,\log(\frac{x_{n-1}}{x_n}))$. The inverse transformation is given by $X=h^{-1}(Y)=(\frac{e^{y_1}}{1+\sum_{i=1}^{n-1}e^{y_i}},\ldots,\frac{e^{y_{n-1}}}{1+\sum_{i=1}^{n-1}e^{y_i}},\frac{1}{1+\sum_{i=1}^{n-1}e^{y_i}})$. As we discussed above, here we will ignore the last component of $X$ when calculate Jacobian of $h^{-1}$. Denote $k=1+\sum_{i=1}^{n-1}e^{y_i}$, the result is 
\begin{align*}
    J_{h^{-1}}&=\frac{1}{k^2}\begin{pmatrix}
        e^{y_1}(k-e^{y_1})&-e^{y_1}e^{y_2}&\ldots&-e^{y_1}e^{y_{n-1}}\\
        -e^{y_2}e^{y_1}&e^{y_2}(k-e^{y_2})&\ldots&\vdots\\
        \vdots&\ldots&\ddots&\vdots\\
        -e^{y_{n-1}}e^{y_1}&-e^{y_{n-1}}e^{y_{2}}&\ldots&e^{y_{n-1}}(k-e^{y_{n-1}})
    \end{pmatrix}\\
    &=\frac{1}{k^2}(k\cdot \mathrm{diag}(e^{y_1},\ldots,e^{y_{n-1}})-\begin{pmatrix}
        e^{y_1}\\
        \vdots\\
        e^{y_{n-1}}
    \end{pmatrix}
    \begin{pmatrix}
        e^{y_1}&\ldots&e^{y_{n-1}}
    \end{pmatrix})
\end{align*}
The form of $J_{h^{-1}}$ is the same as $J_{g^{-1}}$, thus we can immediately apply Lemma1,
\begin{align*}
    \mathrm{det}(J_{h^{-1}})&=\mathrm{det}(\frac{1}{k}\mathrm{diag}(e^{y_1},\ldots,e^{y^{n-1}}))(1-\frac{1}{k}\begin{pmatrix}
        e^{y_1}&\ldots&e^{y_{n-1}}
    \end{pmatrix}\mathrm{diag}(e^{-y_1},\ldots,e^{-y_{n-1}})\begin{pmatrix}
        e^{y_1}\\
        \vdots\\
        e^{y_{n-1}}
    \end{pmatrix}
)\\
&=\frac{1}{k^{n-1}}\cdot e^{\sum_{i=1}^{n-1}y_i}(1-\frac{1}{k}\begin{pmatrix}
        e^{y_1}&\ldots&e^{y_{n-1}}
    \end{pmatrix}\mathrm{diag}(e^{-y_1},\ldots,e^{-y_{n-1}})\begin{pmatrix}
        e^{y_1}\\
        \vdots\\
        e^{y_{n-1}}
    \end{pmatrix})\\
    &=\frac{1}{k^n}\prod_{i=1}^{n-1}e^{y_i}
\end{align*}
so the PDF of $Y$ is 
\begin{align*}
    f_{Y}(y)&=f_{x}(h^{-1}(y))\cdot |J_{h^{-1}}|\\
    &=\frac{1}{B(\alpha)}\prod_{i=1}^{n-1}(e^{y_i})^{\alpha_{i}}\cdot k^{-\sum_{i=1}^{n}\alpha_i}.
\end{align*}

\section{Normalized Negative Binomial Distribution}
In this section, we present a detailed derivation that the normalized negative binomial distributions conform Beta Binomial distribution. From Townes work we know that the Poisson–Gamma mixture model yields a marginal distribution for the counts that is negative binomial\cite{mosimann1962}. Assume $X$ conforms Poisson distribution with fixed $\Lambda$, that is $X|\Lambda\mathrel{\sim}\mathrm{Poisson}(\Lambda)$ and $\Lambda$ is itself a random variable with gamma distribution $\Lambda\mathrel{\sim}\mathrm{Gamma}(r,\theta)$. Then the marginal distribution of $X\mathrel{\sim}\mathrm{NB}(r,p)$, where $p=\frac{\theta}{1+\theta}$. Now suppose that conditional on $\lambda =(\lambda_1,\ldots,\lambda_n)$, the counts $X_1,\ldots,X_n$ are mutually independent with $X_i|\lambda_i\mathrel{\sim}\mathrm{Poisson}(\lambda_i)$, and $\lambda_1,\ldots,\lambda_n$ are independent with $\lambda_i\mathrel{\sim}\mathrm{Gamma}(r_i,\theta)$, sharing the same scale $\theta$. By the common-scale 
summation of gamma distribution, we know $\lambda'=\sum_{i=1}^{n}\lambda_i\mathrel{\sim}\mathrm{Gamma}(R,\theta)$, $R=\sum_{i=1}^nr_i$. And the Poisson distribution has the similar property that the sum of independent Poisson random variables is also Poisson distribution, that is, $S=\sum_{i=1}^{n}x_i|\lambda\mathrel{\sim}\mathrm{Poisson}(\lambda')$. Therefore, we obtain that the marginal distribution of $S$ is Negative Binomial with parameters $R$ and $p=\frac{\theta}{1+\theta}$.

Next, we will derive that for fixed $\lambda=(\lambda_1,\ldots,\lambda_n)$ and $S=m$, $X=(X_1,\ldots,X_n)$ will be multinomial distribution. Firstly we give the joint PMF of $X$ with fixed $\lambda$,
$$\mathrm{Pr}(x_1,\ldots,x_n|\lambda)=\prod_{i=1}^{n}\frac{e^{-\lambda_i}\cdot\lambda_i^{x_i}}{x_i!}.$$
And we know the PMF of $S=\sum_{i-1}^nx_i$ conditional on $\lambda$ is 
$$\mathrm{Pr}(S=m|\lambda)=\frac{e^{-\lambda'}\cdot(\lambda')^m}{m!}.$$
Then we can deduce PMF of $X$ given $S=m$ and $\lambda$
$$\mathrm{Pr}(X|S=m,\lambda)=\frac{\mathrm{Pr}(X,S=m|\lambda)}{\mathrm{Pr}(S=m|\lambda)}.$$
Suppose $A:=\{X=x\}$, $B:=\{S=m\}$, if $\sum_{i=1}^nx_i=m$, then $A\subseteq B$ and $\mathrm{Pr}(A\cap B|\lambda)=\mathrm{Pr}(A|\lambda)$. While if $\sum_{i=1}^nx_i\neq m$, then $A\cap B=\emptyset$ and $\mathrm{Pr}(A\cap B|\lambda)=0$. Thus $\mathrm{Pr}(A\cap B|\lambda)=\mathrm{Pr}(A|\lambda)$, that is $\mathrm{Pr}(X,S=m|\lambda)=\mathrm{Pr}(X|\lambda)$. Then we have
\begin{align*}
    \mathrm{Pr}(X|S=m,\lambda)&=\frac{\mathrm{Pr}(X|\lambda)}{\mathrm{Pr}(S=m|\lambda)}\\
    &=\frac{\prod_{i=1}^ne^{-\lambda_i}\lambda_i^{x_i}/x_i!}{e^{-\lambda'}(\lambda')^m/m!}\\
    &=\frac{\prod_{i=1}^n\lambda_i^{x_i}/x_i!}{(\lambda')^m/m!}\\
    &=\frac{m!}{\prod_{i=1}^nx_i!}\prod_{i=1}^n(\frac{\lambda_i}{\lambda'})^{x_i},
\end{align*}
which is multinomial distribution.

Since $\lambda_i\overset{\text{ind.}}{\sim}\mathrm{Gamma}(r_i,\theta)$, we have $\pi=(\pi_1,\ldots,\pi_n)=(\frac{\lambda_1}{\lambda'},\ldots,\frac{\lambda_n}{\lambda'})\mathrel{\sim}\mathrm{Dirichlet}(r_1,\ldots,r_n)$. Then we can show that $X|S=m$ conforms Dirichlet-Multinomial distribution. Denote $W=\sum_{i=1}^{n}\lambda_i$, we know $W\mathrel{\sim}\mathrm{Gamma}(R,\theta)$ and $\pi\perp W$\cite{mosimann1962}. Since $\lambda=\pi\cdot W$, we can obtain that $X|S=m, \pi, W$ is the same as $X|S=m,\lambda$, which is multinomial distribution. Moreover, 
\begin{align*}
    \mathrm{Pr}(x|S=m)&=\iint\mathrm{Pr}(x,\pi,w|S=m)dwd\pi\\
    &=\iint\mathrm{Pr}(x|S=m,\pi,W)\cdot f(\pi,w|S=m)dwd\pi.
\end{align*}
Here we need some further steps to deal with $f(\pi,w|S=m)$, firstly we introduce an useful result,
\begin{proposition}
    $\pi$ is independent of $S$.
\end{proposition}
\begin{proof}
    Since $S=m|(\pi,W)\mathrel{\sim}\mathrm{Poisson}(\sum_{i=1}^n\lambda_i)=\mathrm{Poisson}(W)$, we can obtain $S\perp\pi |W$. Then 
    \begin{align*}
        f_{S,\pi}(m,\pi)&=\int f_{S,\pi,W}(m,\pi,w)dw\\
        &=\int\mathrm{Pr}(S=m|\pi,W)\cdot f(\pi,w)dw\\
        &=\int\mathrm{Pr}(S=m|\pi,W)\cdot f(\pi)\cdot f(w)dw\\
        &=f(\pi)\int\mathrm{Pr}(S=m|\pi,W)\cdot f(w)dw\\
        &=f(\pi)\cdot f_S(m),
    \end{align*}
    thus $\pi$ is independent of $S$.
\end{proof}
Now by independence, we can rewrite $f(\pi,w|S=m)$ as 
\begin{align*}
    f(\pi,w|S=m)&=\frac{f(\pi,w,S=m)}{\mathrm{Pr}(S=m)}\\
    &=\frac{f(w,S=m)\cdot f(\pi)}{\mathrm{Pr}(S=m)}\\
    &=\frac{\mathrm{Pr}(S=m|W)\cdot f(w)\cdot f(\pi)}{\mathrm{Pr}(S=m)}.
\end{align*}
And we know that $\mathrm{Pr}(x|S=m,\pi,W)\mathrel{\sim}\mathrm{Multinomial}(m,\pi)$ is independent of $W$, so
\begin{align*}
    \mathrm{Pr}(x|S=m)&=\iint\mathrm{Pr}(x|S=m,\pi,W)\cdot\frac{\mathrm{Pr}(S=m|W)\cdot f(w)\cdot f(\pi)}{\mathrm{Pr}(S=m)}dwd\pi\\
    &=\int_{\Delta_{n-1}}\mathrm{Pr}(x|S=m,\pi,W)\cdot f(\pi)d\pi\cdot\underbrace{\frac{\int f(w)\cdot\mathrm{Pr}(S=m|W)dw}{\mathrm{Pr}(S=m)}}_{=1}\\
    &=\int_{\Delta_{n-1}}\underbrace{\frac{m!}{\prod_{i=1}^nx_i!}\prod_{i=1}^n\pi_i^{x_i}}_{\mathrm{Multinomial}(m,\pi)}\cdot\underbrace{\frac{\Gamma(R)}{\prod_{i=1}^n\Gamma(r_i)}\prod_{i=1}^n\pi_i^{r_i-1}}_{\mathrm{Dirichlet(r)}}d\pi.
\end{align*}
By Dirichlet integral $\int_{\Delta_{n-1}}\prod_{i=1}^nz_i^{\alpha_i-1}dz=\frac{\prod_{i=1}^n\Gamma(\alpha_i)}{\Gamma(\sum_{i=1}^{n}\alpha_i)}$, the above integral can be reformulated as 
\begin{align*}
    \mathrm{Pr}(x|S=m)&=\frac{m!\Gamma(R)}{\prod_{i=1}^nx_i!\prod_{i=1}^n\Gamma(r_i)}\int_{\Delta_{n-1}}\prod_{i=1}^n\pi_i^{x_i+r_i-1}d\pi\\
    &=\frac{m!\Gamma(R)}{\prod_{i=1}^nx_i!\prod_{i=1}^n\Gamma(r_i)}\cdot\frac{\prod_{i=1}^n\Gamma(x_i+r_i)}{\Gamma(m+R)}\\
    &=\frac{m!}{\prod_{i=1}^nx_i!}\cdot\frac{\Gamma(R)}{\Gamma(m+R)}\prod_{i=1}^n\frac{\Gamma(x_i+r_i)}{\Gamma(r_i)},
\end{align*}
which is Dirichlet-Multinomial distribution.

Now if we fix one component, and merge the remaining $n-1$ components into a single category, we can then get the marginal distribution $\mathrm{Pr}(x_1=k|S=m)$. To start derivation, we recall that $X|S=m,\pi\mathrel{\sim}\mathrm{Multinomial}(m,\pi)$ and $\pi\mathrel{\sim}\mathrm{Dirichlet}(r)$. We also need to mention that the marginal distribution of Dirichlet is Beta distribution, thus $f_{\pi_1}(u)=\frac{u^{r_1-1}(1-u)^{R-r_1-1}}{B(r_1,R-r_1)}$. Next we fix $x_1=k$ and sum all other components, that is 
\begin{align*}
    \mathrm{Pr}(x_1=k|S=m,\pi_1=u)&=\sum_{\sum_{i=2}^nx_i=m-k}\frac{m!}{k!\prod_{i=2}^nx_i!}\cdot u^k\prod_{i=2}^n\pi_i^{x_i}\\
    &=\frac{m!}{k!}\cdot u^k\sum_{\sum_{i=2}^nx_i=m-k}\frac{1}{\prod_{i=2}^nx_i!}\prod_{i=2}^n\pi_i^{x_i}.
\end{align*}
By Multinomial Theorem, we have $$(\pi_2+\ldots+\pi_n)^{m-k}=\sum_{\sum_{i=2}^nx_i=m-k}\frac{(m-k)!}{\prod_{i=2}^nx_i!}\prod_{i=2}^n\pi_i^{x_i},$$
so $\sum_{\sum_{i=2}^nx_i=m-k}\frac{1}{\prod_{i=2}^nx_i!}\prod_{i=2}^n\pi_i^{x_i}=\frac{(\pi_2+\ldots+\pi_n)^{m-k}}{(m-k)!}=\frac{(1-u)^{m-k}}{(m-k)!}$. Substitute back, we can obtain
\begin{align*}
    \mathrm{Pr}(x_1=k|S=m,\pi_1=u)&=\frac{m!}{k!}\cdot u^k\cdot\frac{(1-u)^{m-k}}{(m-k)!}\\
    &=\begin{pmatrix}
        m\\
        k
    \end{pmatrix}u^k(1-u)^{m-k}
\end{align*}
is Binomial Distribution. Then
$$\mathrm{Pr}(x_1=k|S=m)=\int_0^1\mathrm{Pr}(x_1=k|S=m,\pi_1=u)\cdot f_{\pi_1|S}(u|m)du.$$
By Proposition 1 we know $\pi$ is independent of $S$, thus $f_{\pi_1|S}(u|m)=f_{\pi_1}(u)$. Then the above is equal to 
$$\mathrm{Pr}(x_1=k|S=m)=\int_0^1\begin{pmatrix}
    m\\
    k
\end{pmatrix}u^k(1-u)^{m-k}\cdot\frac{u^{r_1-1}(1-u)^{R-r_1-1}}{B(r_1,R-r_1)}du.$$
Notice that the Beta function $B(z_1,z_2)=\int_0^1t^{z_1-1}(1-t)^{z_2-1}dt$ appears in this integral, it follows that we can further express it as
\begin{align*}
    \mathrm{Pr}(x_1=k|S=m)&=\begin{pmatrix}
        m\\
        k
    \end{pmatrix}\frac{1}{B(r_1,R-r_1)}\int_0^1u^{k+r_1-1}\cdot (1-u)^{R+m-k-r_1-1}du\\
    &=\begin{pmatrix}
        m\\
        k
    \end{pmatrix}\frac{B(k+r_1,R+m-k-r_1)}{B(r_1,R-r_1)}
\end{align*}
called Beta-Binomial distribution.

Finally, take $Y_1=\frac{X_1}{S}$, $y_1\in\{\frac{k}{m},m\geq 1, 0\leq k\leq m\}$. We have the PMF of $Y_1$,
\begin{align*}
    \mathrm{Pr}(Y_1=\frac{k}{m})&=\underbrace{\mathrm{Pr}(S=m)}_{\mathrm{NB}(R,p)}\cdot\underbrace{\mathrm{Pr}(x_1=k|S=m)}_{\mathrm{Beta-Binomial}(m;r_1,R-r_1)}\\
    &=\begin{pmatrix}
        m+R-1\\
        m
    \end{pmatrix}(1-p)^R\cdot p^m\cdot\begin{pmatrix}
        m\\
        k
    \end{pmatrix}\frac{B(k+r_1,m-k+R-r_1)}{B(r_1,R-r_1)}.
\end{align*}

\section{Discussion}
In this article, we assume the scale of Gamma distributions is the same. But in practice, this situation is unusual. A more general assumption is Gamma distributions with integer shape parameter and distinct scale\cite{mathai1981}. In this case, the sum of $n$ gamma random variables will conform a linear combination of $K$ Gamma distributions. Then the marginal distribution $S$ will be the same combination of Negative Binomial. While $X$ conditional on $S$ and $\lambda$ is still Multinomial distribution. Integrate the latent parameter $\pi$ that is now Liouville-Dirichlet, we can obtain $X$ is Liouville-Dirichlet Multinomial under fixed total count $S$. But the marginal is hard to calculate, which need further work. 

\end{document}